\documentclass{elsarticle}

\usepackage{lineno,hyperref}
\modulolinenumbers[5]

\journal{Journal of Differential Equations}

\usepackage{latexsym}
\usepackage{amsmath}
\usepackage{amssymb}
\usepackage{amsfonts}
\usepackage{verbatim}
\usepackage[latin1]{inputenc}
\usepackage{mathrsfs}
\usepackage{amsfonts}
\usepackage{graphicx}
\usepackage{colortbl,dcolumn}
\usepackage{amsmath}
\usepackage{psfrag}
\usepackage{booktabs}

\newtheorem{tm}{Theorem}[section]
\newtheorem{rk}{Remark}[section]

\newtheorem{prop}{Proposition}[section]
\newtheorem{lm}{Lemma}[section]
\newtheorem{cor}{Corollary}[section]

\newcommand{\cc}{\mathbb C}
\newcommand{\ee}{\mathbb E}
\newcommand{\hh}{\mathbb H}
\newcommand{\pp}{\mathbb P}
\newcommand{\nn}{\mathbb N}
\newcommand{\rr}{\mathbb R}
\newcommand{\zz}{\mathbb Z}

\newcommand{\bi}{\mathbf i}
\newcommand{\bs}{\mathbf s}

\newcommand{\CC}{\mathcal C}
\newcommand{\LL}{\mathcal L}

\newcommand{\FFF}{\mathscr F}
\newcommand{\OOO} {\mathscr O}
\newcommand{\<}{\langle}
\renewcommand{\>}{\rangle}

\begin{document}

\begin{frontmatter}

\title{Strong Convergence Rate of Finite Difference Approximations for Stochastic Cubic Schr\"odinger Equations}
\tnotetext[mytitlenote]{This work was supported by National Natural Science Foundation of China (No. 91630312, No. 91530118 and No. 11290142).}

\author[cas]{Jianbo Cui}
\ead{jianbocui@lsec.cc.ac.cn}

\author[cas]{Jialin Hong}
\ead{hjl@lsec.cc.ac.cn}

\author[cas]{Zhihui Liu\corref{cor}}
\ead{liuzhihui@lsec.cc.ac.cn}

\cortext[cor]{Corresponding author.}

\address[cas]{State Key Laboratory of Scientific and Engineering Computing, Institute of Computational Mathematics and Scientific/Engineering Computing, Academy of Mathematics and Systems Science, Chinese Academy of Sciences, Beijing 100190, People's Republic of China}

\begin{abstract}
In this paper, we derive a strong convergence rate of spatial finite difference approximations for both focusing and defocusing stochastic cubic Schr\"odinger equations driven by a multiplicative $Q$-Wiener process. 
Beyond the uniform boundedness of moments for high order derivatives of the exact solution, the key requirement of our approach is the exponential integrability of both the exact and numerical solutions.
By constructing and analyzing a Lyapunov functional and its discrete correspondence, we derive the uniform boundedness of moments for high order derivatives of the exact solution and the first order derivative of the numerical solution, which immediately yields the well-posedness of both the continuous and discrete problems.
The latter exponential integrability is obtained through a variant of a criterion given by [Cox, Hutzenthaler and Jentzen, arXiv:1309.5595].
As a by-product of this exponential integrability, we prove that the exact and numerical solutions depend continuously on the initial data and obtain a 
large deviation-type result on the dependence of the noise with first order strong convergence rate.
\end{abstract}

\begin{keyword}
stochastic cubic Schr\"odinger equation \sep 
strong convergence rate \sep 
central difference scheme \sep 
exponential integrability \sep
continuous dependence
\MSC[2010] 
60H35 \sep 
60H15 \sep 
60G05
\end{keyword}

\end{frontmatter}

\linenumbers


\section{Introduction and main idea}

There is a general theory on strong error estimations for stochastic partial differential equations (SPDEs) with Lipschitz coefficients (see e.g. \cite{ACLW16, BJK16, CHL16, CHL15a} and references therein), where one usually adopts the semigroup or equivalent Green's function framework.
For SPDEs with non-Lipschitz but monotone coefficients, one can use the variational framework to derive strong convergence rates of numerical approximations (see e.g. \cite{KLL15}).
Unfortunately, the monotonicity assumption is too restrictive in the sense that the coefficients of the majority of nonlinear SPDEs from applications, including stochastic Navier-Stokes equations, stochastic Burgers equations, Cahn-Hilliard-Cook equations and stochastic nonlinear Schr\"odinger equations etc., do not satisfy the monotonicity assumption.
Since these SPDEs can not be solved explicitly, one needs to develop effective numerical techniques to study them.
Depending on the particular physical model, it may be necessary to design numerical schemes for solutions of the underlying SPDEs with strong convergence rates, i.e., rates in $\mathbb L^p(\Omega)$ for some $p\ge 1$.
We point out that strong convergence rates are particularly important for efficient multilevel Monte Carlo methods (see e.g. \cite{BFRS16, Keb05}).

Our main purpose in this paper is to derive the strong convergence rate of a representative, spatial approximation for the one-dimensional stochastic cubic Schr\"odinger equation
\begin{align*}
\bi du+(\Delta u+\lambda |u|^2 u) dt &=u\circ dW(t)
\quad \text{in}\quad  (0,T]\times \OOO
\end{align*}
with the initial datum $u(0)=u_0$ under homogenous Dirichlet boundary condition.
Here $\lambda=1$ or $-1$ corresponds to focusing or defocusing cases, $T>0$ is a fixed real number, $\OOO=(0, 1)$ and $W=\{W(t):\ t\in [0,T]\}$ is a $Q$-Wiener process on a stochastic basis $(\Omega, \FFF, \FFF_t, \pp)$, i.e., there exists an real-valued, orthonormal basis $\{e_k\}_{k=1}^\infty $ of $\mathbb L^2(\OOO; \rr)$ and a sequence of mutually independent, real-valued Brownian motions $\{\beta_k\}_{k=1}^\infty $ such that
$W(t)=\sum_{k=1}^\infty Q^{\frac12}e_k\beta_k(t)$ for $t\in [0,T]$.
For convenience, we always consider the equivalent It\^o equation
\begin{align}\label{nls}
du=\Big(\bi \Delta u+\bi \lambda |u|^2 u-\frac12 F_Q u\Big)dt-\bi udW(t)
\quad \text{in}\quad (0,T]\times \OOO 
\end{align}
with the initial datum $u(0)=u_0$,
where $F_Q:=\sum\limits_{k=1}^\infty  (Q^\frac12 e_k)^2$.

The stochastic cubic Schr\"odinger equation \eqref{nls} has been studied in e.g. \cite{DDM02, RGYBC95} to motivate the possible role of noise to prevent or delay collapse formation. 
Many authors concern numerical approximations for Eq. \eqref{nls}; see e.g. \cite{CHP16, BD04, DD06} and references therein. 
To deal with the nonlinearity, one usually apply the truncation technique.
However, this method only produces convergence rates in certain sense such as in probability or pathwise  which is weaker than strong sense. 
In this paper, by deducing the uniform boundedness for moments of high order derivatives of the exact solution and exponential moments of both the exact and numerical solutions, we obtain the strong error estimate of the central difference scheme (see \eqref{dif-nls} below). 
To the best of our knowledge, this is the first result deriving strong convergence rates of numerical approximations for Eq. \eqref{nls}.
For other types of SPDEs, we are only aware that \cite{Dor12} 
analyzes
the spectral Galerkin approximation for the two-dimensional stochastic Navier-Stokes equation using a uniform bound for exponential moments of the solution derived in \cite{HM06}, 
and that \cite{HJ14} 
studies
spectral Galerkin approximations for the one-dimensional Cahn-Hillard-Cook equation and stochastic Burgers equation through exponential integrability characterized in \cite{CHJ13}.

Before proposing our main idea, we introduce some frequently used notations and related technical analytic tools.

\subsection{Notations and analytic tools}

\begin{enumerate}
\item
We denote $\zz_N=\{0,1,\cdots,N\}$ for a given $N\in \nn_+$.
Let $0=x_0<x_1<\cdots<x_{N+1}=1$ be the uniform partition of the interval $\OOO$ with the step size $h=\frac1{N+1}$.
For a grid function $f$, we denote $f(l):=f(x_l)$ for simplicity.
We use $\delta_+$ and $\delta_-$ to denote the forward difference operator and backward difference operator, respectively, i.e., 
\begin{align*}
\delta_+f(l):=\frac {f(l+1)-f(l)}h,\quad 
\delta_-f(l):=\frac {f(l)-f(l-1)}h.
\end{align*}

\item
Denote by $\mathbb L^2:=\mathbb L^2(\OOO; \cc)$ and $l^2_h$ the continuous and discrete Hilbert spaces with inner products
\begin{align*}
\<f, g\>_{\mathbb L^2}:
=\Re \Bigg[\int_{\OOO} \overline f(x) g(x) dx\bigg], \quad 
\<f, g\>_h:=\sum_{l=0}^{N+1} \Re \Big[\overline f(l) g(l)\Big] h.
\end{align*}
We also use the discrete $l^4_h$- and $l^\infty_h$-spaces with norms
\begin{align*}
\|f\|_{l^4_h}:=\bigg(\sum_{l=0}^{N+1}|f(l)|^4h\bigg)^\frac14, \quad 
\|f\|_{l^\infty_h}:=\sup_{l\in \zz_{N+1}} |f(l)|.
\end{align*}

\item
Denote by $\LL_2(\mathbb L^2, \widetilde{H})$ the space of Hilbert-Schmidt operators  $\Psi$ from $\mathbb L^2$ to another separable Hilbert space $\widetilde{H}$, endowed with the norm
\begin{align*} 
\|\Psi\|_{\LL_2(\mathbb L^2, \widetilde{H})}
:=\bigg(\text{tr}\Big[\Psi^*\Psi \Big]\bigg)^\frac12
=\bigg(\sum_{k=1}^\infty  \|\Psi e_k\|_{\widetilde{H}}^2\bigg)^\frac12.
\end{align*} 
In particular, $\LL_2(\mathbb L^2, \hh^\bs)$ is denoted simply by $\LL_2^\bs$.
Here and what follows, $\bs$ is an positive integer and $\hh^\bs$ denotes the usual Sobolev space $H^\bs$, which consists of functions $f$ such that $\nabla^k f$ exist and are square integrable for all $k\in \zz_\bs$, with the additional boundary condtion
$\Delta^j v= 0$ on $\partial \OOO$ for $\nn \ni j<\bs/2$.

\item
To bound the $\|\cdot\|_{\mathbb L^\infty}$-norm, we need the Gagliardo-Nirenberg inequality \begin{align}\label{gn}
\|f\|_{\mathbb L^\infty}^2
&\le 2\|f\|_{\mathbb L^2} \cdot\|\nabla f\|_{\mathbb L^2},
\quad f\in \hh^1
\end{align}
and its discrete correspondence
\begin{align}\label{gn-dis}
\|f^h\|_{l_h^\infty}^2
&\le 2\|f\|_h \|\delta_+ f\|_h,
\quad f\in l^2_h,\ f(0)=f(N+1)=0,\ \delta_+ f\in l^2_h. 
\end{align}
We also use the Sobolev embedding 
$H^1 \hookrightarrow \mathbb L^\infty$:  
\begin{align}\label{sob}
\|f\|_{\mathbb L^\infty}
\le C_0 \|f\|_{H^1}, \quad f\in H^1
\end{align}
for some constant $C_0$ and the fact that $\hh^\bs$ is an algebra, i.e., for any $f,g\in \hh^\bs$, there exists a constant $C=C(\bs)$ such that 
\begin{align}\label{sob-alg}
\|fg\|_{\hh^\bs}\le C\|f\|_{\hh^\bs}\|g\|_{\hh^\bs}.
\end{align}
Throughout the paper $C$ is a generic constant, independent of the discretization parameter $h$, which
will be different
from line to line.
\end{enumerate}

\subsection{Main idea}

Our main aim is to derive the convergence of the spatial central difference scheme 
\begin{align}\label{dif-nls}
du^h(l)
&=\Big(\bi \delta_+ \delta_-u^h(l)+\bi\lambda |u^h(l)|^2 u^h(l) 
-\frac12 F_Q(l) u^h(l)\Big) dt-\bi u^h(l) dW(t, l)
\end{align}
towards Eq. \eqref{nls} with an algebraic rate in strong sense.
Define $R^h(l):
=\Delta u(l)-\delta_+ \delta_-u(l)$ for $l\in \zz_{N+1}$.
The initial datum of Eq. \eqref{dif-nls} is the grid function 
$u^h_0=\{u_0(l),\ l\in \zz_{N+1}\}$.
It is clear that Eq. \eqref{nls} and  Eq. \eqref{dif-nls} possess the continuous and discrete charge conservation laws, respectively, i.e., for all $t\in [0,T]$ it holds a.s. that 
\begin{align}
\label{ccl}
\|u(t)\|^2_{\mathbb L^2}
=\|u_0\|_{\mathbb L^2}^2,\quad
\|u^h(t)\|^2_h
=\|u_0^h\|^2_h.
\end{align} 
The exact solution of Eq. \eqref{nls}, at the grid points, satisfies
\begin{align*}
du(l)
&=\Big(\bi \delta_+ \delta_-u(l)+\bi R^h(l) +\bi \lambda|u(l)|^2 u(l)
-\frac12 u(l)F_Q(l)\Big) dt-\bi u(l) dW(t, l).
\end{align*}
Denote by $\epsilon^h$ the difference between the exact solution $u$ and the numerical solution $u^h$ defined at the grid points $x_l$, $l\in \zz_{N+1}$.
Applying It\^{o} formula to the functional $\|\epsilon^h(t)\|^2_h$, using the continuous and discrete Gagliardo-Nirenberg inequalities \eqref{gn}--\eqref{gn-dis} and charge conservation laws 
\eqref{ccl} (more details see Theorem \ref{u-uh}), we obtain
\begin{align} \label{err-inq}
\|\epsilon^h(t)\|_h^2  
\le \int_0^t \|R^h\|_{l_h^\infty}^2 dr
+\int_0^t \Big(1+2\|u_0\|_{\mathbb L^2}\|\nabla u\|_{\mathbb L^2}
+2\|u_0^h\|_h\|\delta_+ u^h\|_h \Big) \|\epsilon^h\|_h^2 dr.
\end{align}
Taking expectation, as in the deterministic case, leads to  
\begin{align*}
&\ee\bigg[\|\epsilon^h(t)\|_h^2 \bigg] 
\le \int_0^t \ee\bigg[\|R^h\|_{l_h^\infty}^2\bigg] dr \\
&\quad +\int_0^t \ee\bigg[\Big(1+2\|u_0\|_{\mathbb L^2}\|\nabla u\|_{\mathbb L^2}
+2\|u_0^h\|_h\|\delta_+ u^h\|_h \Big) \|\epsilon^h\|_h^2 \bigg] dr.
\end{align*}
Due to the appearance of the nonlinear term in the last integral above, the classical Gronwall inequality is not available and one could not derive a strong convergence rate. 
These difficulties are common features to obtain strong convergence rates for numerical approximations appearing in many situations, see e.g. \cite{CHP16, DD06} for stochastic nonlinear Schr\"odinger equations and \cite{BJ13, KLM11} for other SPDEs with non-monotone coefficients.

Our main idea is applying Gronwall-Bellman inequality to \eqref{err-inq} before taking expectation.
Then using H\"older and Minkowski inequalities, we have
\begin{align}\label{str-emh}
\bigg(\ee\bigg[\sup_{t\in [0,T]}\|\epsilon^h(t)\|_h^2\bigg]\bigg)^\frac12
&\le C\bigg(\ee\bigg[\sup_{t\in [0,T]} \|R^h\|_{l_h^\infty}^4\bigg]\bigg)^\frac14 \nonumber \\ 
&\qquad \bigg\|\exp\bigg(\int_0^T \|u_0\|_{\mathbb L^2}\|\nabla u\|_{\mathbb L^2} dr\bigg) \bigg\|_{\mathbb L^8(\Omega)} \nonumber \\
&\qquad\quad \bigg\|\exp\bigg(\int_0^T \|u_0^h\|_h\|\delta_+ u^h\|_h dr \bigg)\bigg\|_{\mathbb L^8(\Omega)}.
\end{align}
In order to obtain the strong convergence rate for our scheme \eqref{dif-nls}, we need to estimate the three terms appearing on the right hand side of  \eqref{str-emh}.
The first expectation produces the strong convergence rate $\mathcal O(h^2)$ with a multiple given by $p$-moments of the solution under $\hh^\bs$-norm, which is proved to be uniformly bounded in Theorem \ref{us} and Corollary \ref{us-sup}.
To control the last two exponential moments with the random initial datum $u_0$, we apply a variant of a criterion given by \cite{CHJ13} on exponential integrability of a Hilbert-valued stochastic process which is the strong solution of a stochastic differential equation in Hilbert space (see Lemma \ref{exp-int-f0} and Theorem \ref{exp-u-uh}).

Meanwhile, we derive the continuous dependence, in $\mathbb L^p(\Omega; \CC([0,T]; \mathbb L^2))$ with $p=2$ or $p\ge 4$, for the solution of the stochastic cubic Schr\"odinger equation \eqref{nls} on both the initial data and the noises with explicit rates (see Corollaries \ref{con-dep-ini} and \ref{con-dep-noi}).
Similar continuous dependence on the initial data for the numerical solution of the central difference scheme \eqref{dif-nls} can also be obtained.
Such continuous dependence property is not a trivial property for the solutions of SPDEs with non-Lipschitz coefficients. 
We also illustrate that this continuous dependence property is very crucial, besides for theoretical analysis such as in \cite[Chapter 9.1]{DZ14}, for numerical computation because there exist round-off errors in computer simulations. 

The rest of the paper is organized as follows. 
In Section \ref{sec-wel-reg}, we bound the $p$-moments for high order derivatives of the exact solution and discrete first derivative of the numerical solution.
The uniform bound on exponential moments of energy functionals of solutions as well as continuous dependence on initial data and noises is proved in Section \ref{sec-exp-mom}.
The results in Section \ref{sec-wel-reg} and Section \ref{sec-exp-mom} are used in Section \ref{sec-str-ord} to derive the strong convergence rate of the central difference scheme \eqref{dif-nls}.

\section{Well-posedness and regularity}
\label{sec-wel-reg}

In this section, we first prove the moments' uniform boundedness for the solution of the stochastic cubic Schr\"odinger equation \eqref{nls} by analyzing the Lyapunov functional defined by \eqref{lya}, which is necessary for proving the global well-posedness of Eq. \eqref{nls}. 
This uniform boundedness is also useful to derive a strong convergence rate of the central difference scheme \eqref{dif-nls}.
Then  we show, with the help of the discrete energy functional defined by \eqref{Uh}, a priori estimate and thus the well-posedness of the discrete equation \eqref{dif-nls}.

\subsection{A priori estimation of the exact solution}

For $\bs=1$ or $2$, it is known that the stochastic cubic Schr\"odinger equation \eqref{nls} possesses a unique mild solution $u$ which is in $\CC([0,T]; \hh^\bs)$ a.s. under some assumptions on $u_0$ and $Q$ (see \cite{BD03} for $\bs=1$, and \cite{CHP16} for $\bs=2$ in the defocusing case).
Our main result in this part is a priori estimation of $u$ in 
$\hh^\bs$-norm with integer $\bs\ge 2$ for both focusing and defocusing cases.
This will enables us to bound the term $\ee\big[ \|R^h\|_{l_h^\infty}^4\big]$ appearing in \eqref{str-emh}.
We remark that our arguments can also be applied to the whole line.

Throughout this section, we assume that the initial datum $u_0$ is $\FFF_0$-measurable and belongs to $\hh^\bs$ a.s. for certain $\bs\ge 2$. 
To control the nonlinear term $|u|^2 u$ in Eq. \eqref{nls}, we introduce the Lyapunov functional 
\begin{align}\label{lya}
f(u)=\|\nabla^\bs u\|_{\mathbb L^2}^2
-\lambda \big\<(-\Delta)^{\bs-1} u, |u|^2 u \big\>_{\mathbb L^2},
\quad u\in \hh^\bs.
\end{align}
By the inequality \eqref{sob-alg} and integration by parts,
we have 
\begin{align}\label{f0}
\bigg|\ee \bigg[f(u_0)\bigg]\bigg|
&\le C\bigg( \ee\bigg[|u_0|_{\hh^\bs}^2\bigg]
+\ee\bigg[ \|u_0\|^4_{\hh^{\bs-1}}\bigg] \bigg).
\end{align}
Simple calculations yield that the first and the second order derivatives of $f(u)$ are
\begin{align}
Df(u)(v)
&=2 \big\<\nabla^\bs u, \nabla^\bs v \big\>_{\mathbb L^2}
-2\lambda \big\<(-\Delta)^{\bs-1}u, u \Re\big[\overline uv\big] \big\>_{\mathbb L^2} \nonumber  \\
&\quad -\lambda \big\<(-\Delta)^{\bs-1} u, |u|^2 v \big\>_{\mathbb L^2}
-\lambda \big\<(-\Delta)^{\bs-1} v, |u|^2 u \big\>_{\mathbb L^2}, \label{df} \\
D^2f(u)(v,w)
&=2 \big\<\nabla^\bs  v, \nabla^\bs w \big\>_{\mathbb L^2}
-2\lambda \big\<(-\Delta)^{\bs-1} u, w\Re\big[\overline uv\big]\big\>_{\mathbb L^2}
 \nonumber  \\
&\quad-2\lambda \big\<(-\Delta)^{\bs-1} w, u\Re\big[\overline uv\big]\big\>_{\mathbb L^2} 
-2\lambda \big\<(-\Delta)^{\bs-1} u,u\Re\big[\overline vw\big]\big\>_{\mathbb L^2}  
 \nonumber \\
&\quad -2\lambda \big\<(-\Delta)^{\bs-1} u, v\Re\big[\overline uw\big]\big\>_{\mathbb L^2}
-\lambda \big\<(-\Delta)^{\bs-1} w, |u|^2v \big\>_{\mathbb L^2}  \nonumber \\
&\quad -2\lambda \big\<(-\Delta)^{\bs-1} v,u\Re\big[\overline uw\big]\big\>_{\mathbb L^2}
-\lambda \big\<(-\Delta)^{\bs-1} v, |u|^2 w \big\>_{\mathbb L^2}, \label{d2f}
\end{align}
where $v,w\in \CC_0^\infty$.
Applying  the It\^o formula to the functional $f(u)$ defined by \eqref{lya}, we get 
\begin{align}\label{f-f0}
f(u(t))-f(u_0) 
&=\int_0^tDf(u)\bigg(\bi \Delta u+\bi \lambda|u|^2u
-\frac12 F_Q u\bigg)dr  \nonumber  \\
&\quad +\frac12 \int_0^t\text{tr}\bigg[ D^2f(u) \Big (-\bi uQ^{\frac12} \Big)  \Big(-\bi uQ^{\frac12} \Big)^* \bigg] dr  \nonumber \\
&\quad +\int_0^tDf(u)(-\bi u)dW(r)   
=:I_1(t)+I_2(t)+I_3(t).
\end{align}
Substituting the expressions \eqref{df} and \eqref{d2f} of $Df$ and $D^2 f$ into $I_1(t)$ and $I_2(t)$, respectively, we obtain 
\begin{align*}
I_1(t)&=2\int_0^t \Big\<\nabla^\bs  u, \nabla^\bs \Big(\bi\Delta u+\bi\lambda|u|^2u-\frac12 F_Q u\Big) \Big\>_{\mathbb L^2} dr\\
&\quad -\lambda \int_0^t \Big\<(-\Delta)^{\bs-1} u, 
u \Big(\bi\overline u\Delta u-\bi u \Delta\overline u- |u|^2F_Q \Big) \Big\>_{\mathbb L^2} dr\\
&\quad -\lambda \int_0^t \Big\<(-\Delta)^{\bs-1} u, 
|u|^2 \Big(\bi\Delta u+\bi\lambda |u|^2u-\frac12 F_Q u \Big)\Big\>_{\mathbb L^2} dr\\
&\quad -\lambda \int_0^t \Big\<(-\Delta)^{\bs-1}\Big(\bi\Delta u+\bi\lambda |u|^2 u-\frac12 F_Q u\Big), |u|^2u \Big\>_{\mathbb L^2} dr,
\end{align*} 
and
\begin{align*}
I_2(t)
&=\int_0^t \text{tr}\bigg[\nabla^\bs \Big({-\bi uQ^{\frac 12}}\Big)^*\otimes \nabla^\bs \Big(-\bi uQ^{\frac 12} \Big) \bigg]dr\\
&\quad -2\lambda\int_0^t \text{tr}\bigg[\Big((-\Delta)^{\bs-1}\overline u\Big) u\Re \Big[ \Big({-\bi uQ^{\frac 12}} \Big)^*\otimes  \Big(-\bi uQ^{\frac 12} \Big) \Big] \bigg]dr\\
&\quad -\lambda\int_0^t \text{tr}\bigg[(-\Delta)^{\bs-1} \Big( {-\bi uQ^{\frac 12}} \Big)^* \otimes|u|^2  \Big({-\bi uQ^{\frac 12}} \Big) \bigg]dr.
\end{align*}

Our main result in this subsection is the following a priori estimation of algebraic moments for high order derivatives of the solution of Eq. \eqref{nls}.

\begin{tm}\label{us}
Let $p=2$ or $p\ge 4$  and $\bs \ge 2$. 
Assume that 
\begin{align}\label{u0}
u_0\in \bigcap\limits_{m=2}^\bs \mathbb L^{3^{\bs-m}p}(\Omega; \hh^m)\cap \bigcap\limits_{m=0}^1 \mathbb L^{3^{\bs-m-1}5p}(\Omega; \hh^m)
\end{align} 
and $Q^\frac12 \in \LL_2^{ \bs}$. 
Then there exists a constant $C=C(T,p,u_0,Q)$ such that 
\begin{align}\label{uhs}
\sup_{t\in[0,T]}\ee\bigg[ \|u(t)\|_{\hh^\bs}^p\bigg]
\le C.
\end{align}
\end{tm}

\textbf{Proof}
Let $t\in [0,T]$.
We first estimate $I_1(t)$ in \eqref{f-f0}.
Denote the four integrals in $I_1(t)$ successively by $I_{11}(t)$, $I_{12}(t)$, 
$I_{13}(t)$ and $I_{14}(t)$.
Integration by parts yields that $I_{11}(t)$ and $I_{12}(t)$ can be rewritten as 
\begin{align*}
I_{11}(t)= I_{111}(t)+I_{112}(t), \quad
I_{12}(t)= I_{121}(t)+I_{122}(t),
\end{align*} 
where 
\begin{align*}
I_{111}(t)
&=-2\int_0^t \Big\<(-\Delta)^{\bs-1} u, \bi \lambda \Delta \Big(|u|^2u\Big)\Big\>_{\mathbb L^2} dr, \\
I_{112}(t)
&=-\int_0^t \Big\<\nabla^\bs  u, 
\nabla^\bs \Big(F_Q u\Big) \Big\>_{\mathbb L^2} dr, \\
I_{121}(t)
&=-\lambda\int_0^t \Big\<(-\Delta)^{\bs-1} u,
u \Big(\bi\overline u\Delta u-\bi  u\Delta \overline u \Big) \Big\>_{\mathbb L^2} dr,  \\
I_{122}(t)
&=\lambda\int_0^t \Big\< (-\Delta)^{\bs-1} u, |u|^2 F_Q u \Big\>_{\mathbb L^2} dr.
\end{align*} 
The Cauchy-Schwarz inequality and the inequality \eqref{sob-alg} yield that
\begin{align*}
\bigg|\ee\bigg[I_{112}(t)\bigg]\bigg|
+\bigg|\ee\bigg[I_{122}(t)\bigg]\bigg| 
&\le C \bigg(\int_0^t \ee\bigg[\|u(r)\|_{\hh^{\bs-1}}^4\bigg] dr
+\int_0^t \ee\bigg[|u(r)|_{\hh^\bs}^2\bigg] dr\bigg).
\end{align*}
The term $I_{111}(t)$ is divided into two equal parts which can balance 
$I_{121}(t)$ and $I_{14}(t)$.
More precisely, inserting the identities 
$\Delta (|u|^2 u) =2 |u|^2 \Delta u+4u |\nabla u|^2
+2\overline u (\nabla u)^2+u^2\Delta \overline u$ 
and $\<\Delta^{\bs-1}(|u|^2u), \bi \Delta u\>_{\mathbb L^2}
+\<\Delta^{\bs-1}u, \bi \Delta (|u|^2u)\>_{\mathbb L^2}=0$, we obtain
\begin{align*}
& I_{111}(t)+I_{121}(t)+I_{14}(t) \\
&=\bigg(\frac{I_{111}(t)}2+I_{121}(t)\bigg)
+\bigg(\frac{I_{111}(t)}2+I_{14}(t)\bigg) \\
&=-\lambda \int_0^t \Big\<(-\Delta)^{\bs-1} u, 
3\bi |u|^2\Delta u \Big\>_{\mathbb L^2} dr \\
&\quad -\lambda \int_0^t \bigg\<(-\Delta)^{\bs-1} u, 
4\bi |\nabla u|^2u+2\bi (\nabla u)^2\overline u \bigg\>_{\mathbb L^2} dr  \\
&\quad -\lambda\int_0^t \bigg\<(-\Delta)^{\bs-1}\Big(|u|^2u\Big),
-\frac12 F_Q u \bigg\>_{\mathbb L^2} dr 
=:I_a(t)+I_b(t)+I_c(t).
\end{align*}
Applying integration by parts, Leibniz formula and the fact that $\<\nabla^{\bs}u, \bi |u|^2\nabla^\bs u\>_{\mathbb L^2}=0$, we have 
\begin{align*}
\Big\<(-\Delta)^{\bs-1}u, \bi |u|^2 \Delta u \Big\>_{\mathbb L^2} 
&=\Big\<\nabla^{\bs}u, 
\bi \nabla^{\bs-2} \Big(|u|^2 \Delta u\Big) \Big\>_{\mathbb L^2} \\
&=\sum_{j=0}^{\bs-3}C_{\bs-2}^j 
\Big\<\nabla^{\bs}u, \nabla^{\bs-2-j}\Big(|u|^2\Big) \cdot \nabla^j\Big(\Delta u\Big) \Big\>_{\mathbb L^2}.
\end{align*}
Then by the Sobolev embedding \eqref{sob} and the inequality \eqref{sob-alg} we get for $\bs>2$,
\begin{align*}
\Big|\Big\<(-\Delta)^{\bs-1}u, \bi |u|^2 \Delta u\Big\>_{\mathbb L^2}\Big| 
\le C|u|_{\hh^{\bs}} \sum_{j=0}^{\bs-3} 
\|\nabla^{\bs-2-j}|u|^2\|_{\mathbb L^\infty} \|u\|_{\hh^{j+2}}  
\le C|u|_{\hh^{\bs}} \|u\|_{\hh^{\bs-1}}^3.
\end{align*}
The above estimate is also valid for $\bs=2$, since 
$\<\Delta u, \bi |u|^2 \Delta u\>_{\mathbb L^2}=0$.
This implies that
\begin{align*}
\bigg|\ee\bigg[I_a(t)\bigg]\bigg|
\le C \bigg(\int_0^t \ee\bigg[\|u(r)\|_{\hh^{\bs-1}}^6\bigg] dr
+\int_0^t \ee \bigg[|u(r)|^2_{\hh^{\bs}}\bigg] dr \bigg).
\end{align*}
Applying H\"older inequality, integration by parts and the inequality \eqref{sob-alg}, we obtain 
for $\bs>2$,
\begin{align*}
\bigg|\ee\bigg[I_b(t)\bigg]\bigg|
&\le C \bigg(\int_0^t \ee\bigg[\|u(r)\|_{\hh^{\bs-1}}^6\bigg] dr
+\int_0^t \ee \bigg[|u(r)|^2_{\hh^{\bs}}\bigg] dr \bigg).
\end{align*} 
When $\bs=2$, by using the Sobolev embedding \eqref{sob}, the Gagliardo-Nirenberg inequality \eqref{gn} and 
the Young inequality, we get
\begin{align*}
\bigg|\ee\bigg[I_b(t)\bigg]\bigg|
&\le \int_0^t \ee\bigg[ |u|_{\hh^2} \|u\|_{\mathbb L^\infty}
\|\nabla u\|_{\mathbb L^\infty}\|\nabla u\|_{\mathbb L^2}\bigg] dr \\
&\le C  \bigg(\int_0^t \ee\bigg[\|u(r)\|_{\hh^{\bs-1}}^{10}\bigg] dr
+\int_0^t \ee \bigg[|u(r)|^2_{\hh^{\bs}}\bigg] dr \bigg).
\end{align*}
Similar arguments imply that
\begin{align*}
\bigg|\ee\bigg[I_c(t)\bigg]\bigg|
&\le C \int_0^t \ee\bigg[\|u(r)\|_{\hh^{\bs-1}}^4\bigg] dr.
\end{align*}
As a result, there exists a constant $C=C(T,Q)$ such that   
\begin{align*}
&\bigg|\ee\bigg[ I_{111}(t)+I_{121}(t)+I_{14}(t)\bigg]\bigg|\\
&\le C  \bigg(1+\int_0^t \ee\bigg[\|u(r)\|_{\hh^{\bs-1}}^{10}\bigg] dr
+\int_0^t \ee \bigg[|u(r)|^2_{\hh^{\bs}}\bigg] dr \bigg).
\end{align*}
For $I_{13}(t)$, using integration by parts and the inequality \eqref{sob-alg}, we have
\begin{align*}
 \bigg|\ee \bigg[I_{13}(t) \bigg] \bigg|
&\le C \bigg(\ee \bigg[|I_a(t) | \bigg]
+\int_0^t\ee \bigg[ \bigg|\Big\<(-\Delta)^{\bs-1}u,  \bi |u|^4 u+\frac12 F_Q |u|^2 u \Big\>_{\mathbb L^2} \bigg| \bigg]dr  \bigg) \\
&\le C \bigg(1+\int_0^t \ee\bigg[\|u(r)\|_{\hh^{\bs-1}}^6\bigg] dr
+\int_0^t \ee \bigg[|u(r)|^2_{\hh^{\bs}}\bigg] dr\bigg).
\end{align*}
Combining the estimations of $I_{11}$ to $I_{14}$, we have 
\begin{align}\label{i1-1}
 \bigg|\ee \bigg[I_1(t) \bigg] \bigg|
&\le C \bigg(1+\int_0^t \ee\bigg[\|u(r)\|_{\hh^{\bs-1}}^6\bigg] dr
+\int_0^t \ee \bigg[|u(r)|^2_{\hh^{\bs}}\bigg] dr\bigg)
\end{align}
for $\bs>2$ and 
\begin{align}\label{i1-2}
 \bigg|\ee \bigg[I_1(t) \bigg] \bigg|
&\le C \bigg(1+\int_0^t \ee\bigg[\|u(r)\|_{\hh^1}^{10}\bigg] dr
+\int_0^t \ee \bigg[|u(r)|^2_{\hh^2}\bigg] dr\bigg).
\end{align}

Now we turn to the estimations of $I_2(t)$ and $I_3(t)$ in \eqref{f-f0}.
Using the Cauchy-Schwarz inequality and the inequality \eqref{sob-alg}, we get  
\begin{align}\label{i2}
 \bigg|\ee \bigg[I_2(t) \bigg] \bigg|
&\le C \bigg(1+\int_0^t \ee\bigg[\|u(r)\|_{\hh^{\bs-1}}^6\bigg] dr
+\int_0^t \ee \bigg[|u(r)|^2_{\hh^{\bs}}\bigg] dr\bigg).
\end{align}
On the other hand, owing to the property of It\^o integral, we have
\begin{align}\label{i3}
 \Big|\ee  \Big[I_3(t)  \Big]  \Big|
&=0.
\end{align}

For $p\ge 4$, we apply the It\^o formula to $f^{\frac p2}(u)$ and obtain
\begin{align}\label{fp}
f^{\frac p2}(u(t))
&=f^{\frac p2}(u_0)+\frac p2\int_0^t  f^{{\frac p2}-1}(u) Df(u)
\Big(\bi \Delta u+\bi \lambda|u|^2u-\frac12 F_Q u \Big) dr \nonumber \\
&\quad + \frac {p(p-2)} 8\int_0^t
\text{tr}\bigg[   f^{{\frac p2}-2}(u) Df(u)
\Big(-\bi uQ^{\frac12} \Big) Df(u)\Big(-\bi uQ^{\frac12}\Big )^*\bigg]dr \nonumber \\
&\quad+\frac p4\int_0^t \text{tr}\bigg[  f^{{\frac p2}-1}(u) D^2f(u)
\Big(-\bi uQ^{\frac12} ,-\bi uQ^{\frac12}\Big) \bigg]dr
\nonumber\\
&\quad+\frac p2 \int_0^t f^{{\frac p2}-1}(u) Df(u) \Big(-\bi u \Big) dW(r).
\end{align}
It follows from the inequality \eqref{sob-alg} and the Cauchy-Schwarz inequality that
\begin{align*}
f^{{\frac p2}-1}(u)
\le C \Big(\|u\|_{\hh^{\bs-1}}^{2(p-2)}+|u|_{\hh^\bs}^{p-2} \Big).
\end{align*}
By the estimations \eqref{i1-1} and \eqref{i1-2} of $I_1(t)$, it holds a.s. that  
\begin{align*}
\Big|Df(u)\Big(\bi \Delta u+\bi \lambda|u|^2u-\frac12 F_Q u\Big)\Big| 
&\le C\Big(\|u\|_{\hh^{\bs-1}}^6+|u|_{\hh^\bs}^2\Big)
\end{align*}
for $\bs>2$ and 
\begin{align*}
\Big|Df(u)\Big(\bi \Delta u+\bi \lambda|u|^2u-\frac12 F_Q u\Big)\Big| 
&\le C\Big(\|u\|_{\hh^1}^{10}+|u|_{\hh^2}^2\Big).
\end{align*}
Then by the Young inequality, we get an estimate for the first integral in \eqref{fp}:
\begin{align*}
& \ee\bigg[\bigg|\int_0^t  f^{{\frac p2}-1}(u) Df(u)
\Big(\bi \Delta u+\bi \lambda|u|^2u-\frac12 F_Q u \Big)dr\bigg| \bigg]\\
&\le C \int_0^t \bigg(\ee\bigg[\|u\|_{\hh^{\bs-1}}^{2p+2}\bigg]
+\ee\bigg[ \|u\|_{\hh^{\bs-1}}^{2p-4} |u|_{\hh^\bs}^2 \bigg]
+\ee\bigg[ \|u\|_{\hh^{\bs-1}}^6 |u|_{\hh^\bs}^{p-2}\bigg]
+\ee \bigg[|u|_{\hh^\bs}^p \bigg] \bigg) dr  \\
&\le C\bigg(1+\int_0^t \ee\bigg[\|u(r)\|_{\hh^{\bs-1}}^{3p}\bigg] dr
+\int_0^t \ee \bigg[|u(r)|_{\hh^\bs}^p \bigg]dr\bigg)
\end{align*}
for $\bs>2$ and
\begin{align*}
& \ee\bigg[\bigg|\int_0^t  f^{{\frac p2}-1}(u) Df(u)
\Big(\bi \Delta u+\bi \lambda|u|^2u-\frac12 F_Q u \Big)dr\bigg| \bigg]\\
&\le C\bigg(1+\int_0^t \ee\bigg[\|u(r)\|_{\hh^1}^{5p}\bigg] dr
+\int_0^t \ee \bigg[|u(r)|_{\hh^2}^p \bigg]dr\bigg).
\end{align*}
Similar arguments can be applied to other terms in \eqref{fp}.
Thus we obtain
\begin{align*}
\ee\bigg[|u(t)|_{\hh^\bs}^p\bigg]
&\le C \bigg(\ee\bigg[f^\frac p2(u(t)) \bigg]
+\ee\bigg[\|u\|_{\hh^{\bs-1}}^{2p}\bigg]\bigg) \\
&\le C\bigg(1+\ee\bigg[\|u_0\|_{\hh^\bs}^p \bigg]
+\ee\bigg[\|u_0\|_{\hh^{\bs-1}}^{2p}\bigg] \\
&\qquad +\int_0^t \ee\bigg[ \|u(r)\|_{\hh^{\bs-1}}^{3p} \bigg] dr
+\int_0^T \ee\bigg[\|u(r)\|_{\hh^\bs}^{p}\bigg] dr\bigg)
\end{align*}
for $\bs>2$ and
\begin{align*}
\ee\bigg[|u(t)|_{\hh^2}^p\bigg]
&\le C\bigg(1+\ee\bigg[\|u_0\|_{\hh^2}^p \bigg]
+\ee\bigg[\|u_0\|_{\hh^1}^{2p}\bigg] \\
&\qquad +\int_0^t \ee\bigg[ \|u(r)\|_{\hh^1}^{5p} \bigg] dr
+\int_0^T \ee\bigg[\|u(r)\|_{\hh^2}^{p}\bigg] dr\bigg).
\end{align*}
Gronwall inequality then yields that
\begin{align*}
\ee\bigg[|u(t)|_{\hh^\bs}^p\bigg]
&\le C\bigg(1+\ee\bigg[\|u_0\|_{\hh^\bs}^p \bigg]
+\ee\bigg[\|u_0\|_{\hh^{\bs-1}}^{2p}\bigg] 
+\int_0^T \ee\bigg[ \|u(r)\|_{\hh^{\bs-1}}^{3p} \bigg] dr\bigg)
\end{align*}
for $\bs>2$ and that 
\begin{align*}
\ee\bigg[|u(t)|_{\hh^2}^p\bigg]
&\le C\bigg(1+\ee\bigg[\|u_0\|_{\hh^2}^p \bigg]
+\ee\bigg[\|u_0\|_{\hh^1}^{2p}\bigg] 
+\int_0^T \ee\bigg[ \|u(r)\|_{\hh^1}^{5p} \bigg] dr\bigg).
\end{align*}
Similar arguments as in \cite[Theorems 4.1 and 4.6]{BD03} lead to 
\begin{align*}
\sup_{t\in [0,T]}\ee\bigg[ \|u(t)\|_{\hh^1}^{5p}\bigg]<\infty
\end{align*} 
provided that $u_0\in \mathbb L^{5p}(\Omega; \hh^1)\cap \mathbb L^{15p}(\Omega; \mathbb L^2)$.
This implies that 
\begin{align*}
\sup_{t\in [0,T]}\ee\bigg[ \|u(t)\|_{\hh^2}^p\bigg]<\infty
\end{align*} 
provided that $u_0\in \mathbb L^p(\Omega; \hh^2)\cap \mathbb L^{5p}(\Omega; \hh^1)\cap \mathbb L^{15p}(\Omega; \mathbb L^2)$.
For $\bs=3$, when $u_0\in \mathbb L^p(\Omega; \hh^3)\cap \mathbb L^{3p}(\Omega; \hh^2)\cap \mathbb L^{15p}(\Omega; \hh^1)\cap \mathbb L^{45p}(\Omega; \mathbb L^2)$, it holds that
\begin{align*}
\ee\bigg[|u(t)|_{\hh^3}^p\bigg]
&\le C\bigg(1+\ee\bigg[\|u_0\|_{\hh^3}^p \bigg]
+\ee\bigg[\|u_0\|_{\hh^2}^{2p}\bigg] 
+\int_0^T \ee\bigg[ \|u(r)\|_{\hh^2}^{3p} \bigg] dr\bigg)
<\infty.
\end{align*}
By induction, we complete the proof of \eqref{uhs}.
\qed\\

\begin{rk}
In \cite[Theorem 4.6]{BD03}, a uniform bound for the Hamiltonian
\begin{align}\label{U}
U(X):=\frac 12\|\nabla X\|_{\mathbb L^2}^2-\frac {\lambda}4 \|X\|_{\mathbb L^4}^4
\end{align}
is used to construct a unique solution with continuous $\hh^1$-valued paths for Eq. \eqref{nls}. 
We can follow the same strategy as in \cite{BD03} to construct the unique local mild solution with continuous $\hh^\bs$-valued paths by a contraction argument, and then show that it is global by the a priori estimate \eqref{uhs} with $\bs\ge 2$. 
To prove that this mild solution is also a strong one of Eq. \eqref{nls}, we refer to \cite{LR15}, Propositions G.0.4 and G.0.5.
\end{rk}

\begin{cor}\label{us-sup}
Under the same conditions of Theorem \ref{us}, there exists a constant $C=C(T,p,u_0,Q)$ such that 
\begin{align}\label{uhs-sup}
\ee\bigg[\sup_{t\in[0,T]} \|u(t)\|_{\hh^\bs}^{p} \bigg]\le C.
\end{align}
\end{cor}

\textbf{Proof.}
The main step to derive \eqref{uhs-sup} from \eqref{uhs} is that we need to deal with the stochastic integral $I_3(t)$ which is vanished in Theorem \ref{us}.
By the expression \eqref{df} of $Df(u)$, we get 
\begin{align*}
I_3(t)
&=2\int_0^t \Big\<\nabla^\bs u, \nabla^\bs \Big(-\bi udW(r)\Big) \Big\>_{\mathbb L^2} \\
&\quad -\lambda\int_0^t \Big\<(-\Delta)^{\bs-1} u, |u|^2 \Big(-\bi udW(r)\Big) \Big\>_{\mathbb L^2} \\
&\quad -\lambda\int_0^t \Big\<(-\Delta)^{\bs-1}\Big(-\bi u dW(r)\Big), 
|u|^2u \Big\>_{\mathbb L^2}.
\end{align*}
Applying the Burkholder-Davis-Gundy, H\"older and Young inequalities and using similar arguments to estimate $I_a(t)$ in Theorem \ref{us}, we obtain
\begin{align*}
& \ee\bigg[\sup_{t\in[0,T]} |I_3(t)|^{\frac p2}\bigg]  \\
&\le 
C \bigg( \ee\bigg[\bigg|\int_0^T |u(r)|_{\hh^\bs}^2 \|u(r)\|_{\hh^{\bs-1}}^2 dr \bigg|^{\frac p4}\bigg]
+\ee\bigg[\bigg| \int_0^T |u(r)|_{\hh^\bs}^2 \|u(r)\|_{\hh^{\bs-2}}^6 dr \bigg|^{\frac p4}\bigg]\bigg)\\
&\le C\ee \Bigg[\bigg(\sup_{t\in[0,T]}|u(t)|_{\hh^\bs}^{\frac p2}\bigg)
\cdot \bigg (\int_0^T\|u(t)\|_{\hh^{\bs-1}}^6 dt \bigg)^{\frac p4} \Bigg]\\
&\le \frac 1{2^{p+1}} \ee \bigg[\sup_{t\in[0,T]}|u(t)|_{\hh^\bs}^p\bigg]
+C\int_0^T\ee \bigg[\|u(t)\|_{\hh^{\bs-1}}^{3p}\bigg]dt.
\end{align*}
Combining the above estimate and the estimations of 
\begin{align*}
\ee\bigg[\sup_{t\in[0,T]} |I_1(t)|^{\frac p2}\bigg]
\quad \text{and}\quad 
\ee\bigg[\sup_{t\in[0,T]} |I_2(t)|^{\frac p2}\bigg]
\end{align*} 
derived in Theorem \ref{us}, we conclude the estimate \eqref{uhs-sup} .
\qed\\

\subsection{A priori estimation of the numerical solution}

The local existence and uniqueness of the  solution for the central difference scheme \eqref{dif-nls} can be proved by the contraction argument in \cite{BD03} for Eq. \eqref{nls}. 
Then global existence is an immediate consequence of the following a priori estimate.
To this end, we define the discrete energy functional 
\begin{align}\label{Uh}
U^h(u^h):
=\frac 12 \|\delta_+u^h\|_{h}^2-\frac {\lambda}4\|u^h\|_{l^4_h}^4
\end{align}
similarly to the continuous one $U$ defined by \eqref{U}.

\begin{prop}\label{mom-bou-uh}
Let $p=2$ or $p\ge 4$.
Assume that $u_0^h \in \mathbb L^{3p}(\Omega, l^2_h)$, $\delta_+ u_0^h \in \mathbb L^p(\Omega, l^2_h)$ and $Q\in \LL_2^2$. 
Then there exists a constant $C=C(T,p,u_0,Q)$ such that 
\begin{align}\label{h1h-sup}
\ee\bigg[\sup_{t\in[0,T]}\|\delta_+u^h(t)\|_h^{p}\bigg]
\le C.
\end{align}
\end{prop}

\textbf{Proof.}
We only prove \eqref{h1h-sup} for $p=2$ since the proof for the case $p\ge 4$ is similar to that of Theorem \ref{us} and Corollary \ref{us-sup}.
Applying the It\^o formula to the energy functional $U^h(u^h)$ defined by \eqref{Uh}, we obtain
\begin{align*}
& U^h(u^h(t))-U^h(u_0)  \\
&=\lambda \int_0^t \Big\<\delta_+u^h, \bi \delta_+\Big(|u^h|^2 u^h\Big) \Big\>_h dr 
-\sum_{k=1}^\infty  \int_0^t
\Big\<\delta_+u^h, \bi u^h \delta_+ \Big(Q^{\frac 12}e_k\Big) \Big\>_h d\beta_k(r)   \\
&\quad +\frac12 \Big\<\delta_+ \delta_- u^h, u^h F_Q\Big\>_h
+\lambda \int_0^t \Big\<\delta_+ \delta_- u^h, \bi|u^h|^2 u^h \Big\>_h dr \\
&\quad +\frac 12\int_0^t\sum_{k=1}^\infty 
\Big\|\delta_+ \Big(u^h Q^{\frac 12}e_k\Big) \Big\|_h^2dr 
:=U_a +U_b+U_c+U_d+U_e.
\end{align*}
Due to the symmetry of the numerical scheme \eqref{dif-nls} under the Dirichlet boundary condition, the term $U_a+U_d$ vanishes.
Then the martingale property of the It\^o integral yields that
$\ee[ U_b]=0$ and thus 
\begin{align*}
\ee\Big[ U^h(u^h(t))\Big]
=\ee\Big[U^h(u_0)\Big]+\ee\Big[U_c+U_e\Big].
\end{align*}
Applying integration by parts, we obtain
\begin{align*}
& U_c+U_e \\
&=-\frac 12\int_0^t \bigg(\sum_{l=0}^{N+1} h\Re\Big[\overline{\delta_+u^h(l)}\delta_+u^h(l)F_Q(l+1)\Big]+
\Big\<\delta_+u^h(l), u^h(l) \delta_+F_Q(l) \Big\>_h \bigg)dr\\
&\quad +\sum_{k=1}^\infty \frac 12\int_0^t
\bigg(\sum_{l=0}^{N+1} h|\delta_+u^h(l)Q^{\frac12}e_k(l+1)|^2
+\Big\|u^h \delta_+\Big(Q^{\frac 12}e_k\Big)(l)\Big\|_h^2\bigg)dr\\
&\quad +\sum_{k=1}^\infty\int_0^t \sum_{l=0}^{N+1}h\Re\Big[\overline{\delta_+u^h(l)}u^h(l)\delta_+\Big(Q^{\frac12}e_k\Big)(l)Q^{\frac12}e_k(l+1)\Big]dr.
\end{align*}
Similar calculations as in Theorem \ref{us} deduce that 
\begin{align*}
\ee\Big[ U^h(u^h(t)) \Big]
&\le \ee\Big[U^h(u_0)\Big]
  +\frac {3t}2 \sum_{k=1}^\infty  \|\nabla (Q^\frac12 e_k)\|_{\mathbb L^\infty}^2\ee\Big[\|u^h_0\|_h^2\Big] \\
 &\le \ee\Big[U^h(u_0)\Big]
+  \frac {3C_0^2t}2 \|Q^\frac12\|_{\LL_2^2}^2
\ee\Big[\|u^h_0\|_h^2\Big],
\end{align*}
where $C_0$ is the Sobolev embedding coefficient in \eqref{sob}.
The Cauchy-Schwarz inequality and the discrete Gagliardo-Nirenberg inequality \eqref{gn-dis}  imply that 
\begin{align}\label{eq}
\frac14\|\delta_+ u^h\|_h^2-\frac14 \|u^h\|_{h}^6
\le U^h(u^h)
\le \frac 34\|\delta_+ u^h\|_h^2+\frac14 \|u^h\|_{h}^6.
\end{align}
As a result, we get 
\begin{align*}
\sup_{t\in[0,T]}\ee\Big[U^h(u^h(t))\Big]
&\le \frac 34 \ee\Big[\|\delta_+u^h_0\|_h^2\Big]
+\frac 14\ee\Big[\|u_0^h\|_{h}^6\Big]
+C\|Q^{\frac 12}\|_{\LL_2^2}^2 \ee\Big[\|u_0^h\|_{h}^2\Big]T.
\end{align*}
By the charge conservation laws \eqref{ccl} and the inequality \eqref{eq}, there exists a constant $C=C(T,Q,C_0)$ such that 
\begin{align*}
\sup_{t\in[0,T]}\ee\Big[\|\delta_+u^h\|^2_h\Big]
&\le C\Big(\ee\Big[\|\delta_+u^h_0\|_h^2\Big]
+\ee\Big[\|u_0^h\|_{h}^6\Big]
+\ee\Big[\|u^h_0\|_{h}^2\Big]\Big).
\end{align*}
We conclude the uniform estimate \eqref{h1h-sup} by similar arguments to derive \eqref{uhs-sup} from \eqref{uhs} as in Corollary \ref{us-sup}.
\qed\\

\begin{rk}
Similar arguments of \cite[Lemma 5 and Lemma 6]{CHP16} yield the existence of continuous versions of both $u$ and $u^h$ under the assumptions in Theorem \ref{us} and Proposition \ref{mom-bou-uh}.
This continuity will be used in the next section to derive the exponential integrability properties of both $u$ and $u^h$.
\end{rk}

\section{Exponential integrability and continuous dependence}
\label{sec-exp-mom}

In this section, we establish the exponential integrability for both the stochastic cubic Schr\"odinger equation \eqref{nls} and its central difference scheme \eqref{dif-nls}.
This exponential integrability is used in the next section to bound the last two exponential moments in \eqref{str-emh}.
As a by-product, the exact solution  depends continuously on the initial data as well as on the noise, under $\mathbb L^p(\Omega; \CC([0,T]; \mathbb L^2))$-norm for any $p=2$ or $p\ge 4$, with explicit rate.

\subsection{Exponential integrability property}

To handle the uniform boundedness of the last two exponential moments in \eqref{str-emh} for the solutions of Eq. \eqref{nls} and Eq. \eqref{dif-nls}, we give a criterion on exponential integrability, which is a variant of \cite[Corollary 2.4]{CHJ13}.

%
%

\begin{lm}\label{exp-int-f0}
Let $H$ be a Hilbert space and $X$ be an $H$-valued adapted stochastic process with continuous sample paths satisfying 
$\int_0^T\|\mu(X_t)\|+\|\sigma(X_t)\|^2 dt<\infty$ a.s., and for all $t\in [0,T]$, 
$X_t=X_0+\int_0^t \mu(X_r)dr +\int_0^t \sigma(X_r)dW_r$ a.s.
Assume that there exist two functionals $\overline V$ and $V\in \CC^2(H;\rr)$ and an $\FFF_0$-measurable random variable $\alpha$ such that for almost every $t\in [0,T]$,
\begin{align}\label{con-exp}
&DV(X_t)\mu(X_t)
+\frac{\text{tr}\big[D^2V(X_t)\sigma(X_t)\sigma^*(X_t)\big]}2 \nonumber  \\
&\quad +\frac{\|\sigma^*(X) DV(X_t)\|^2}{2e^{\alpha t}}+\overline V(X_t)
\le \alpha V(X_t),\quad a.s.
\end{align}
Then 
\begin{align}\label{exp-int}
\sup_{t\in [0,T]}\ee\bigg[\exp\bigg( \frac {V(X_t)}{e^{\alpha t}}+\int_0^t\frac {\overline V(X_r)}{e^{\alpha r}}dr \bigg)\bigg]
\le \ee\bigg[ e^{V(X_0)} \bigg].
\end{align}
\end{lm}

\textbf{Proof.}
Let $Y_t=\int_0^t\frac {\overline V(X_r)}{e^{\alpha r}}dr$.
Applying the It\^o formula to 
\begin{align*}
Z(t,X_t,Y_t):=\exp\bigg(\frac {V(X_t)}{e^{\alpha t}}+Y_t \bigg), 
\end{align*}
we obtain 	
\begin{align*}	
& Z(t,X_t,Y_t)-e^{V(X_0)} \\
&=\int_0^t e^{-\alpha r} Z(r,X_r,Y_r) (\overline V(X_r)-\alpha V(X_r)) dr \\
&\quad +\int_0^t e^{-\alpha r} Z(r,X_r,Y_r) DV(X_r)dX_r\\
&\quad +\frac 12 \int_0^t e^{-\alpha r} Z(r,X_r,Y_r)
\text{tr}\Big[ D^2V(X_r)\sigma(X_r) \sigma^*(X_r)\Big] dr\\
&\quad +\frac 12\int_0^t e^{-2\alpha r} Z(r,X_r,Y_r) 
\big\|\sigma^*(X_r)DV(X_r)\big\|^2dr.
\end{align*}
Condition \eqref{con-exp} implies that 
\begin{align*}
&\ee \Big[Z(t,X_t,Y_t)\Big]-\ee \Big[e^{V(X_0)}\Big] \\
&=\ee\bigg[\int_0^t Z(r,X_r,Y_r)\bigg(\frac{\overline V(X_r)-\alpha V(X_r)}{e^{\alpha r}}+\frac{DV(X_r)\mu(X_r)}{e^{\alpha r}} \\
&\qquad +\frac{\text{tr}\bigg[D^2V(X_r)\sigma(X_r) \sigma^*(X_r)\bigg]}{2e^{2\alpha r}}+\frac{\|\sigma^*(X_r)DV(X_r)\|^2}
{2e^{\alpha r}}\bigg)dr\bigg]
\le 0.
\end{align*} 
This completes the proof of \eqref{exp-int}.
\qed\\

In the rest of this section, we assume that the stochastic cubic Schr\"odinger equation \eqref{nls} and the central difference scheme \eqref{dif-nls} possess unique strong solutions with continuous trajectories. 

Applying Lemma \ref{exp-int-f0} to the energy functionals $U$ and $U^h$ defined by \eqref{U} and \eqref{Uh}, respectively, we obtain the following uniform bounds of exponential moments for $u$ and $u^h$.
This is the main ingredient in Section \ref{sec-str-ord} to deduce the strong error estimate between $u$ and $u^h$.

\begin{tm}\label{exp-u-uh}
Let $q\ge 1$ and $Q^\frac12\in \LL_2^2$.
Assume that 
\begin{align}\label{con-u1}
\ee\bigg[ e^{U(u_0)}\bigg]
+\ee\bigg[\exp\bigg(\frac{\|u_0\|_{\mathbb L^2}^6}2
+4q^2  T^2 \|u_0\|_{\mathbb L^2}^2
e^{2 C_0^2\|Q^\frac12\|_{\LL_2^2}^2 \|u_0\|_{\mathbb L^2}^2T} \bigg)\bigg]<\infty
\end{align}
and 
\begin{align}\label{con-u2}
\ee\bigg[ e^{U^h(u_0^h)}\bigg]
+\ee\bigg[\exp\bigg(\frac{\|u_0^h\|_h^6}2
+4q^2  T^2 \|u_0^h\|_h^2
e^{2 C_0^2\|Q^\frac12\|_{\LL_2^2}^2 \|u_0^h\|_h^2T}\bigg)\bigg]<\infty.
\end{align}
Then there exist a constant $C=C(T,q,u_0,Q)$ such that 
\begin{align}\label{exp-mom-u}
\bigg\|\exp\bigg(\int_0^T \|u_0\|_{\mathbb L^2}\|\nabla u\|_{\mathbb L^2}dr\bigg)\bigg\|_{\mathbb L^q(\Omega)}
&\le C, \\
\label{exp-mom-uh}
\bigg\|\exp\bigg(\int_0^T \|u_0^h\|_h\|\delta_+ u^h\|_h dr\bigg)\bigg\|_{\mathbb L^q(\Omega)}
&\le C.
\end{align}
\end{tm}

\textbf{Proof.}
Simple calculations as \eqref{df} and \eqref{d2f} show that
\begin{align*}
DU(X)Y
&= \<\nabla X, \nabla Y\>_{\mathbb L^2}-\lambda \<|X|^2 X, Y\>_{\mathbb L^2}, \\
D^2 U(X)(Y,Z)
&=\<\nabla Z, \nabla Y\>_{\mathbb L^2}-\lambda \<|X|^2 Y,Z\>_{\mathbb L^2}
-2\lambda \<\Re[\overline XY]X, Z\>_{\mathbb L^2}.
\end{align*}
In the case of Eq. \eqref{nls}, $\mu(u)=\bi \Delta u+\bi\lambda |u|^2u-\frac12F_Q u$ and $\sigma(u)=-\bi uQ^\frac12$.
Then
\begin{align*}
&DU(u)\mu(u)
=-\frac12 \Big\<\nabla u, \nabla F_Q u\Big\>_{\mathbb L^2}
-\frac12 \Big\<\nabla u, u\nabla F_Q \Big\>_{\mathbb L^2}
+\frac\lambda 2\<|u|^4, F_Q\>_{\mathbb L^2}, \\
&\text{tr}\Big[D^2 U(u) \sigma(u)\sigma^*(u)\Big] 
=\sum_{k=1}^\infty \Big\|\nabla \Big(u Q^\frac12 e_k\Big) \Big\|_{\mathbb L^2}^2
-{\lambda}\Big\<|u|^4, F_Q \Big\>_{\mathbb L^2},
\end{align*}
and
\begin{align*}
\Big\|\sigma^*(u)DU(u)\Big\|_{\mathbb L^2}^2 
=\Big\< \nabla u, -\bi u \sum_{k=1}^\infty \nabla \Big(Q^\frac12e_k \Big) \Big\>_{\mathbb L^2} ^2.
\end{align*}
Therefore, by the Sobolev embedding \eqref{sob}, we have 
\begin{align*}
& DU(u)\mu(u)+\frac12 \text{tr}\Big[D^2U(u) \sigma(u)\sigma^*(u)\Big]
+\frac 1{2e^{\alpha t}}\|\sigma^*(u)DU(u)\|^2 \\
&=\frac12 \sum_{k=1}^\infty  
\Big\<|u|^2, \Big(\nabla Q^\frac12 e_k\Big)^2 \Big\>_{\mathbb L^2} 
+\frac 1{2e^{\alpha t}} \sum_{k=1}^\infty  
\Big\< \nabla u, -\bi u \nabla \Big(Q^\frac12e_k\Big) \Big\>_{\mathbb L^2} ^2\\
&\le \frac{C_0^2}2\|Q^\frac12\|_{\LL_2^2}^2 \|u_0\|_{\mathbb L^2}^2 
+\frac{C_0^2}2\|Q^\frac12\|_{\LL_2^2}^2 \|u_0\|_{\mathbb L^2}^2\|\nabla u\|_{\mathbb L^2}^2 .
\end{align*}
We conclude that for $\lambda=-1$,
\begin{align*}
&DU(u)\mu(u)+\frac{\text{tr} [D^2U(u) \sigma(u)\sigma^*(u)]}{2}
+\frac{\|\sigma^*(u)DU(u)\|^2}{2e^{\alpha t}}
\le \alpha_{-1} U(u)+\beta_{-1}
\end{align*}
with $\alpha_{-1}=C_0^2 \|Q^\frac12\|_{\LL_2^2}^2\|u_0\|_{\mathbb L^2}^2$ and 
$\beta_{-1} =\frac{C_0^2}2 \|Q^\frac12\|_{\LL_2^2}^2 \|u_0\|_{\mathbb L^2}^2$. 
When $\lambda=1$, using the fact that 
$\|u\|_{\mathbb L^4}^4\le 2\|u\|_{\mathbb L^2}^3 \|\nabla u\|_{\mathbb L^2}$, we get
\begin{align*}
&DU(u)\mu(u)+\frac{\text{tr}[D^2U(u)\sigma(u)\sigma^*(u)]}{2}
+\frac{\|\sigma^*(u)DU(u)\|^2}{2e^{\alpha t}}
\le \alpha_1 U(u)+\beta_1
\end{align*}
with $\alpha_1=2C_0^2\|Q^\frac12\|_{\LL_2^2}^2 \|u_0\|_{\mathbb L^2}^2$ and 
$\beta_1=\frac {C_0^2}2\|Q^\frac12\|_{\LL_2^2}^2 (\|u_0\|_{\mathbb L^2}^2+\|u_0\|_{\mathbb L^2}^8)$.
Applying Lemma \ref{exp-int-f0} with $\overline{U}=-\beta_\lambda$ for $\lambda=\pm 1$, we obtain 
\begin{align}\label{lem-nls}
\sup_{t\in [0,T]}\ee\bigg[\exp\bigg(\frac {U(u(t))}{e^{\alpha_\lambda t}}-\int_0^t \frac {\beta_\lambda}{e^{\alpha_\lambda s}}ds\bigg)\bigg]
\le \ee\bigg[ e^{U(u_0)}\bigg].
\end{align}
When $\lambda = -1$, \eqref{lem-nls} yields that
\begin{align*}
\sup_{t\in[0,T]}\ee\bigg[\exp\bigg(\frac {U(u(t))}{e^{\alpha_{-1}t}}\bigg)\bigg]
\le  e^{\frac 12}\ee\bigg[ e^{U(u_0)}\bigg].
\end{align*}
Applying the Young inequality and a variant of Jensen inequality, we obtain 
\begin{align*}
&\sup_{t\in [0,T]}\bigg\|\exp\bigg(\int_0^t \|u_0\|_{\mathbb L^2} 
\|\nabla u\|_{\mathbb L^2} ds\bigg) \bigg\|_{\mathbb L^q(\Omega)}  \\
&\le \bigg\|\exp\bigg( \frac T{2\epsilon} \|u_0\|_{\mathbb L^2}^2 \bigg) \bigg\|_{\mathbb L^{2q}(\Omega)}
\bigg\|\exp\bigg(\int_0^T \frac \epsilon 2\|\nabla u\|_{\mathbb L^2} ^2ds\bigg) \bigg\|_{\mathbb L^{2q}(\Omega)}  \\
&\le  \bigg\|\exp\bigg( \frac T{2\epsilon} \|u_0\|_{\mathbb L^2}^2 \bigg) \bigg\|_{\mathbb L^{2q}(\Omega)}
\sup_{t\in[0,T]}\bigg\|\exp\bigg(T\epsilon U(u(t))\bigg)\bigg\|_{\mathbb L^{2q}(\Omega)}.
\end{align*}
Let $\epsilon=\frac 1{2qTe^{\alpha_{-1}T}}$, then the above two estimations yield that
\begin{align*}
&\sup_{t\in [0,T]}\bigg\|\exp\bigg(\int_0^t \|u_0\|_{\mathbb L^2} \|\nabla u\|_{\mathbb L^2} ds\bigg)\bigg\|_{\mathbb L^q(\Omega)} \\
&\le e^{\frac 1{4q}}\sqrt[2q]{\ee\bigg[\exp\bigg({2q^2 T^2 e^{\alpha_{-1} T}\|u_0\|_{\mathbb L^2}^2}\bigg)\bigg]} 
\sqrt[2q]{\ee\bigg[ e^{U(u_0)}\bigg]}.
\end{align*}

When $\lambda=1$, by the fact that
$U(u)\ge \frac14(\|\nabla u\|_{\mathbb L^2}^2-\|u_0\|_{\mathbb L^2}^6)$ and a version of Jensen inequality, we have 
\begin{align*}
&\sup_{t\in [0,T]}\bigg\|\exp\bigg(\int_0^t \|u_0\|_{\mathbb L^2} \|\nabla u\|_{\mathbb L^2} ds\bigg)\bigg\|_{\mathbb L^q(\Omega)}  \\
&\le \bigg\|\exp\bigg(\int_0^T \bigg[\frac{\|u_0\|_{\mathbb L^2}^2}{\epsilon}+\frac{\epsilon\|u_0\|_{\mathbb L^2}^6}{4}
+\frac {\beta_1}{2qe^{\alpha_1r}} \bigg] dr \bigg)\bigg\|_{\mathbb L^{2q}(\Omega)}\\
&\quad \times \bigg\|\exp\bigg(\int_0^T \bigg[
\frac{\epsilon (\|\nabla u\|_{\mathbb L^2}^2-\|u_0\|_{\mathbb L^2}^6)}{4}
-\frac {\beta_1}{2qe^{\alpha_1r}} \bigg] dr \bigg)\bigg\|_{\mathbb L^{2q}(\Omega)}  \\
&\le  \bigg\|\exp\bigg( \frac{T\|u_0\|_{\mathbb L^2}^2}{\epsilon} 
+\frac {\epsilon T \|u_0\|_{\mathbb L^2}^6}{4}
+\frac {\beta_1(1-e^{{-\alpha_1}T})}{2q\alpha_1}\bigg) 
\bigg\|_{\mathbb L^{2q}(\Omega)}\\
&\quad \times \bigg\|\exp\bigg( \epsilon T U(u(t))
-\frac {\beta_1(1-e^{{-\alpha_1}T})}{2q\alpha_1}\bigg) 
\bigg\|_{\mathbb L^{2q}(\Omega)}.
\end{align*}
Let $\epsilon=\frac 1{2qTe^{\alpha_{1}T}}$, then by the Young inequality, 
\begin{align*}
&\sup_{t\in [0,T]}\bigg\|\exp\bigg(\int_0^t \|u_0\|_{\mathbb L^2} \|\nabla u\|_{\mathbb L^2} ds\bigg)\bigg\|_{\mathbb L^q(\Omega)}  \\
&\le \sqrt[2q]{\ee\bigg[\exp\bigg( 4q^2  T^2 e^{\alpha_1T}\|u_0\|_{\mathbb L^2}^2
+\frac{\|u_0\|_{\mathbb L^2}^6}{4e^{\alpha_1T}}
+\frac {\beta_1(1-e^{{-\alpha_1}T})}{\alpha_1} \bigg)\bigg]} \sqrt[2q]{\ee\bigg[e^{U(u_0)}\bigg]} \\
&\le e^{\frac1{8q}}
\sqrt[2q]{\ee\bigg[\exp\bigg(\frac12 \|u_0\|_{\mathbb L^2}^6
+4q^2  T^2 e^{\alpha_1T}\|u_0\|_{\mathbb L^2}^2\bigg)\bigg]} \sqrt[2q]{\ee\bigg[e^{U(u_0)}\bigg]}.
\end{align*}
This shows the estimate \eqref{exp-mom-u}.

For the discrete case, we have
\begin{align*}
D U^h(X^h) Y^h
&=\Big\<\delta_+ X^h, \delta_+ Y^h\Big\>_h
-\lambda \Big\<|X^h|^2 X^h, Y^h \Big\>_h, \\
D^2 U^h(X^h)(Y^h,Z^h)
&=\Big\<\delta_+ Z^h, \delta_+ Y^h \Big\>_h
-\lambda \Big\<|X^h|^2 Y^h,Z^h \Big\>_h\\
&\quad -2\lambda \Big\<\Re\Big[\overline X^hY^h\Big] X^h, Z^h \Big\>_h.
\end{align*}
In the case of Eq. \eqref{dif-nls}, $\mu^h(u^h)(l)=\bi \delta_+ \delta_- u^h(l)+\bi\lambda |u^h(l)|^2u^h(l)-\frac12 F_Q(l) u^h(l)$ and
$\sigma^h(u^h)(e_k)(l)=-\bi u^h(l)Q^{\frac12}e_k(l)$.
Substituting them into the above equalities and using the Sobolev embedding \eqref{sob}, we obtain

\begin{align*}
& DU^h(u^h)\mu^h(u^h)+\frac{\text{tr}[D^2U^h(u^h) \sigma^h(u^h)\sigma^h(u^h)^*]}2
+\frac{\|\sigma^h(u^h)^* DU^h(u^h)\|^2}{2e^{\alpha^h t}} \\
&=-\frac12\sum_{k=1}^{\infty} \Big\<\delta_+ u^h, 
u^h  \Big|\delta_+ \Big(Q^{\frac12}e_k \Big) \Big|^2 h \Big\>_h 
+\frac12\sum_{k=1}^{\infty} \Big\<|u^h|^2, 
 \Big|\delta_+ \Big(Q^{\frac12}e_k \Big) \Big|^2  \Big\>_h \\
&\quad +\frac1{2e^{\alpha^h t}} \sum_{k=1}^{\infty}
 \Big\<\delta_+ u^h, -\bi u^h \delta_+\Big (Q^\frac12e_k \Big) \Big \>_h ^2\\
&\le  \frac {3C_0^2}2 \|Q^\frac12\|_{\LL_2^2}^2
\|u^h_0\|_h^2
+\frac{C_0^2}{2} \|Q^\frac12\|_{\LL_2^2}^2
\|u_0^h\|_h^2\|\delta_+ u^h\|_h^2.
\end{align*}
We conclude that, for $\lambda=-1$,
\begin{align*}
& DU^h(u^h)\mu^h(u^h)+\frac{\text{tr}[D^2U^h(u^h) \sigma^h(u^h)\sigma^h(u^h)^*]}2 \\
&\qquad  +\frac{\|\sigma^h(u^h)^* DU^h(u^h)\|^2}{2e^{\alpha_{-1}^h t}}
\le \alpha_{-1}^h U^h(u^h)+\beta_{-1}^h
\end{align*}
with $\alpha^h_{-1}=C_0^2 \|Q^\frac12\|_{\LL_2^2}^2 \|u_0^h\|_h^2$ and $\beta_{-1}^h=\frac{3C_0^2}{2} \|Q^\frac12\|_{\LL_2^2}^2\|u_0^h\|_h^2$.
Similarly to the continuous case,
\begin{align*}
&\sup_{t\in [0,T]}\bigg\|\exp\bigg(\int_0^t \|u_0^h\|_h \|\delta_+ u^h\|_h ds\bigg)\bigg\|_{\mathbb L^q(\Omega)}  \\
&\le e^{\frac 1{4q}} \sqrt[2q]{\ee\bigg[\exp\bigg(4q^2  T^2
e^{{\alpha^h_{-1}}T} \|u_0^h\|_h^2  \bigg) \bigg]} 
\sqrt[2q]{\ee\bigg[ e^{U^h(u_0^h)}\bigg]}.
\end{align*}
When $\lambda=1$, Similar arguments yield that 
\begin{align*}
& DU^h(u^h)\mu^h(u^h)+\frac{\text{tr}[D^2U^h(u^h) \sigma^h(u^h)\sigma^h(u^h)^*]}2 \\
&\qquad  +\frac{\|\sigma^h(u^h)^* DU^h(u^h)\|^2}{2e^{\alpha_1^h t}}
\le \alpha_1^h U^h(u^h)+\beta_1^h
\end{align*}
with $\alpha^h_1=2C_0^2 \|Q^\frac12\|_{\LL_2^2}^2 \|u_0^h\|_h^2$ and $\beta_{1}^h=\frac {C_0^2}2\|Q^\frac12\|_{\LL_2^2}^2 (3\|u_0\|_{\mathbb L^2}^2+\|u_0\|_{\mathbb L^2}^8)$.
Moreover, we have 
\begin{align*}
&\sup_{t\in [0,T]}\bigg\|\exp\bigg(\int_0^t \|u_0^h\|_h \|\delta_+u^h\|_h ds\bigg)\bigg\|_{\mathbb L^q(\Omega)}  \\
&\le \sqrt[2q]{\ee\bigg[\exp\bigg( 4q^2 T^2 e^{\alpha_1^h T}\|u_0^h\|_h^2
+\frac{\|u_0^h\|_h^6}{4e^{\alpha_1^h T}}
+\frac {\beta^h_1(1-e^{{-\alpha_1^h}T})}{\alpha_1^h} \bigg)\bigg]} \sqrt[2q]{\ee\bigg[e^{U^h(u_0^h)}\bigg]} \\
&\le e^{\frac3{8q}}
\sqrt[2q]{\ee\bigg[\exp\bigg(\frac12 \|u_0^h\|_h^6
+4q^2  T^2 e^{\alpha_1T}\|u_0^h\|_h^2\bigg)\bigg]} \sqrt[2q]{\ee\bigg[e^{U^h(u_0^h)}\bigg]}.
\end{align*}
This shows the estimate \eqref{exp-mom-uh} and we complete the proof.
\qed\\

\subsection{Continuous dependence on initial data and noises}

As an application of the exponential integrability Theorem \ref{exp-u-uh},
we can prove that the solution of the stochastic cubic Schr\"odinger equation \eqref{nls} is continuously depending on the initial data. 
Such continuous dependence property is not a trivial property for the solutions of SPDEs with non-Lipschitz coefficients. 

\begin{cor}\label{con-dep-ini}
Let $p=2$ or $p\ge 4$.
Assume that the condition \eqref{con-u1} holds for $u_0$ and $v_0$ with $q=4p$ and  $Q^\frac12\in \LL_2^2$. 
Let $u=\{u(t):\ t\in [0,T]\}$ and $v=\{v(t):\ t\in [0,T]\}$ be the solutions of Eq. \eqref{nls} with initial data $u_0$ and $v_0$, respectively.
Then there exists a constant $C=C(T,p,u_0,v_0,Q)$ such that
\begin{align}\label{con-dep-ini0}
\bigg(\ee\bigg[\sup_{t\in [0,T]}\|u(t)-v(t)\|_{\mathbb L^2}^p\bigg]\bigg)^\frac1p
\le C \bigg(\ee\bigg[\|u_0-v_0\|_{\mathbb L^2}^{2p}\bigg]\bigg)^\frac1{2p}.
\end{align}
\end{cor}

\textbf{Proof.}
We only prove the case $p=2$ and the proof for the other cases is similar. 
Applying It\^o isometry to the functional $\|u(t)-v(t)\|_{\mathbb L^2}^2$, we obtain
\begin{align*}
&\|u(t)-v(t)\|_{\mathbb L^2}^2-\|u_0-v_0\|_{\mathbb L^2}^2 \\
&=\int_0^t \sum_{k=1}^\infty\|(u-v) Q^\frac12 e_k\|_{\mathbb L^2}^2 dr
+2\int_0^t \Big\<u-v, -\bi (u-v) \Big\>_{\mathbb L^2} dW(s)  \\
&\quad +2\int_0^t \Big\<u-v, \bi \Delta(u-v)+\bi \lambda (|u|^2 u-|v|^2 v)-\frac12 F_Q (u-v) \Big\>_{\mathbb L^2} dr \\
&\le \int_0^t 2 \|u\|_{\mathbb L^\infty} \|v\|_{\mathbb L^\infty}
\|u-v\|_{\mathbb L^2}^2 dr.
\end{align*}
Applying Gronwall inequality and taking expectation, 
combing with the Cauchy-Schwarz inequality and the Gagliardo-Nirenberg inequality \eqref{gn},
we get
\begin{align*}
&\ee\bigg[\sup_{t\in [0,T]}\|u(t)-v(t)\|_{\mathbb L^2}^2\bigg] \\
& \le \bigg\|\exp\bigg(\int_0^T 2\|u\|_{\mathbb L^\infty} \|v\|_{\mathbb L^\infty} dr\bigg) \bigg\|_{\mathbb L^2(\Omega)} 
\bigg(\ee\bigg[\|u_0-v_0\|_{\mathbb L^2}^4\bigg]\bigg)^\frac12\\
&\le \bigg(\ee\bigg[\|u_0-v_0\|_{\mathbb L^2}^4\bigg]\bigg)^\frac12
\bigg\|\exp\bigg(\int_0^T 2\|u\|_{\mathbb L^2} \|\nabla u\|_{\mathbb L^2}\bigg) \bigg\|_{\mathbb L^4(\Omega)} \\
&\qquad \qquad \bigg\|\exp\bigg(\int_0^T 2\|u\|_{\mathbb L^2} \|\nabla u\|_{\mathbb L^2}\bigg) \bigg\|_{\mathbb L^4(\Omega)},
\end{align*}
from which we  conclude \eqref{con-dep-ini0} by Theorem \ref{exp-u-uh}.
\qed\\

\begin{rk}
One can show that the numerical solution $u^h$ of the central difference scheme \eqref{dif-nls} possesses analogous continuous dependence \eqref{con-dep-ini0}, since the essential requirement of exponential integrability of $u^h$ has been established in Theorem \ref{exp-u-uh}.
This continuous dependence property is very crucial for numerical computation, because of round-off errors in computer simulations. 
\end{rk}

Beyond the above continuous dependence result on the initial data, we also have the following 
large deviation-type result, on the stochastic cubic Schr\"odinger equation
\begin{align}\label{eps-nls}
\bi du^\epsilon+\Delta u^\epsilon dt+\lambda |u^\epsilon|^2 u^\epsilon dt=\epsilon u^\epsilon\circ dW,
\quad u^\epsilon(0)=u_0,
\quad \epsilon\in \rr	
\end{align} 
driven by small scale noise.
This type of large deviation estimation can characterize, in $\mathbb L^p(\Omega)$-sense, the deviation of the perturbed solution $u^\epsilon$ from the unperturbed solution $u^0$, which is the solution of Eq. \eqref{eps-nls} with $\epsilon=0$ (see e.g. \cite{FW12}).
 
\begin{cor}\label{con-dep-noi}
Let $p=2$ or $p\ge 4$.
Assume that the condition \eqref{con-u1} holds with $q=4p$ and  $Q^\frac12\in \LL_2^2$. 
Let $u^\epsilon=\{u^\epsilon(t):\ t\in [0,T]\}$ and $u^0=\{u^0(t):\ t\in [0,T]\}$ be the solutions of Eq. \eqref{eps-nls} with the same initial datum $u_0$, respectively.
Then there exists a constant $C=C(T,p,u_0, Q)$ such that
\begin{align}\label{con-dep-noi0}
\bigg(\ee\bigg[\sup_{t\in [0,T]}\|u^\epsilon(t)-u^0(t)\|_{\mathbb L^2}^p\bigg]\bigg)^\frac1p
\le C |\epsilon|.
\end{align}
\end{cor}

\textbf{Proof.}
We only prove the case $p=2$ and similar arguments yield the other cases.
Applying It\^o isometry to the functional $\|u^\epsilon(t)-u^0(t)\|_{\mathbb L^2}^2$, we obtain
\begin{align*}
&\|u^\epsilon(t)-u^0(t)\|_{\mathbb L^2}^2 \\
&=\epsilon^2 \int_0^t \sum_{k}\|u^\epsilon Q^\frac12 e_k\|_{\mathbb L^2}^2 dr
+2\int_0^t \Big\<u^\epsilon-u^0, -\bi (u^\epsilon-u^0) \Big\>_{\mathbb L^2} dW(s)  \\
&\quad +2\int_0^t \Big\<u^\epsilon-u^0, \bi \Delta(u^\epsilon-u^0)+\bi \lambda (|u^\epsilon|^2 u^\epsilon-|u^\epsilon|^2 u^\epsilon)
-\frac12 \epsilon^2 F_Q u^\epsilon  \Big\>_{\mathbb L^2} dr \\
&\le \epsilon^2 \int_0^t \sum_{k} \Big\<u^0, u^\epsilon  |Q^\frac12 e_k|^2 \Big\>_{\mathbb L^2} dr
+\int_0^t 2 \|u^\epsilon\|_{\mathbb L^\infty} \|u^0\|_{\mathbb L^\infty}
\|u^\epsilon-u^0\|_{\mathbb L^2}^2 dr.
\end{align*}
Applying Gronwall-Bellman inequality, taking expectation and using Sobolev embedding, we get 
\begin{align*}
\sqrt{\ee\bigg[\sup_{t\in [0,T]}\|u^{\epsilon}(t)-u^0(t)\|_{\mathbb L^2}^2\bigg]}
\le C \bigg\|\exp\bigg(\int_0^T \|u^\epsilon\|_{\mathbb L^\infty} \|u^0\|_{\mathbb L^\infty} dr\bigg)\bigg\|_{\mathbb L^4(\Omega)}
|\epsilon|.
\end{align*}
Similarly to the proof of Corollary \ref{con-dep-ini}, we conclude \eqref{con-dep-noi0} by  Theorem \ref{exp-u-uh}.
\qed\\

\begin{rk}
For stochastic nonlinear Schr\"odinger equation driven by an additive noise, \cite[Corollary 3.1 and Proposition 3.5]{BD03} derived a.s. continuous dependence, but without convergence rates, on the initial data and the noises. 
\end{rk}

\section{Strong convergence rate of finite difference approximation}
\label{sec-str-ord}

In this section, we use the a priori estimations in Section \ref{sec-wel-reg} and the exponential integrability properties in Section \ref{sec-exp-mom} to establish the strong error estimation between the numerical solution of the spatial central difference scheme 
\begin{align}\label{dif-nls-i}
du^h(l)
&=\Big(\bi \delta_+ \delta_-u^h(l)+\bi\lambda |u^h(l)|^2 u^h(l) 
-\frac12 F_Q(l) u^h(l) \Big) dt -\bi u^h(l) dW(t,l)
\end{align}
and the exact solution of Eq. \eqref{nls}.
It is clear that the exact solution of Eq. \eqref{nls}, at the grid points, satisfies
\begin{align}\label{nls-gri-i}
du(l)
&=\Big(\bi \delta_+ \delta_-u(l)+\bi R^h(l) +\bi \lambda|u(l)|^2 u(l)
-\frac12 F_Q(l) u(l)\Big) dt-\bi u(l) dW(t, l),
\end{align}
where $R^h(l):=\Delta u(l)-\delta_+ \delta_-u(l)$ for $l\in \zz_{N+1}$.

Denote by $\epsilon^h(t,l):=u(t,l)-u^h(t,l)$ the difference between the exact solution $u$ and the numerical solution $u^h$ at the grid point $x_l$ for $t\in [0,T]$ and $l\in \zz_{N+1}$.
Our main result of this paper is the following $\mathbb L^2(\Omega; \CC([0,T]; \mathbb L^2))$-error estimate of $\epsilon^h$.

\begin{tm}\label{u-uh}
Assume that $u_0\in \bigcap\limits_{m=2}^5 \mathbb L^{4\cdot 3^{5-m}}(\Omega; \hh^m)$, \eqref{con-u1}--\eqref{con-u2} hold for  $q=8$ and $Q^\frac12\in \LL_2^{5}$.
Then there exists a constant $C=C(T,u_0,Q)$ such that
\begin{align}\label{u-uh0}
\bigg(\ee\bigg[\sup_{t\in [0,T]}\|u(t)-u^h(t)\|_h^2\bigg]\bigg)^\frac12
&\le C h^2.
\end{align}
\end{tm}

\textbf{Proof.}
Subtracting  Eq. \eqref{dif-nls-i} from  Eq. \eqref{nls-gri-i} and using the identity $|a|^2a-|b|^2b=(|a|^2+|b|^2)(a-b)+ab~ \overline{(a-b)}$
for $a,b\in \mathbb C$, we obtain
\begin{align*}
d\epsilon^h(l)
&=\Big(\bi \delta_+ \delta_- \epsilon^h(l)
+\bi \lambda(|u(l)|^2+|u^h(l)|^2)\epsilon^h(l)
+\bi \lambda u(l)u^h(l)\overline{\epsilon^h}(l) \Big) dt  \\
&\qquad -\frac12 F_Q(l) \epsilon^h(l) dt
+\bi R^h(l) dt
-\bi \epsilon^h(l) dW(t, l).
\end{align*}
Applying It\^{o} formula to the functional $\|\epsilon^h(t)\|^2_h$ and using the fact that $\epsilon^h(0,l)=0$ for any $l\in \zz_{N+1}$ and the symmetry of the proposed scheme under homogeneous Dirichlet boundary condition, we get
\begin{align*}
\sum_{l=0}^{N+1} |\epsilon^h(t, l)|^2
&=2\lambda \int_0^t\sum_{l=0}^{N+1} \Re\bigg[\bi u(l)u^h(l)
\overline{\epsilon^h}(l)^2\bigg]  dr 
+2\int_0^t\sum_{l=0}^{N+1} \Re\bigg[\bi \overline{\epsilon^h}(l) R^h(l)\bigg] dr.
\end{align*}
Applying the Cauchy-Schwarz inequality and the discrete Gagliardo-Nirenberg inequality \eqref{gn-dis}, we obtain
\begin{align*}
\|\epsilon^h(t)\|_h^2  
&\le \int_0^t \|R^h\|_{l_h^\infty}^2 dr
+\int_0^t \Big(1+\|u\|_{\mathbb L^\infty}^2+\|u^h\|_{l_h^\infty}^2 \Big) \|\epsilon^h\|_h^2 dr \\
&\le \int_0^t \|R^h\|_{l_h^\infty}^2 dr
+\int_0^t \Big(1+2\|u_0\|_{\mathbb L^2}\|\nabla u\|_{\mathbb L^2}
+2\|u_0^h\|_h\|\delta_+ u^h\|_h\Big) \|\epsilon^h\|_h^2 dr.
\end{align*}
It follows from Gronwall-Bellman inequality that
\begin{align*}
& \|\epsilon^h(t)\|_h^2 
\le \bigg(\int_0^T \|R^h(r)\|_{l_h^\infty}^2 dr\bigg) \\
&\qquad \times \exp\bigg(\int_0^t \Big[1+2\|u_0\|_{\mathbb L^2}\|\nabla u(r)\|_{\mathbb L^2}+2\|u_0^h(r)\|_h\|\delta_+ u^h(r)\|_h\Big] dr\bigg).
\end{align*}
Taking expectation and using H\"older and Minkovskii inequalities, we derive
\begin{align}\label{err-str}
\bigg(\ee\bigg[\sup_{t\in [0,T]}\|\epsilon^h(t)\|_h^2\bigg]\bigg)^\frac12
&\le T^\frac12 e^\frac T2 \bigg(\ee\bigg[\sup_{t\in [0,T]} \|R^h\|_{l_h^\infty}^4\bigg]\bigg)^\frac14 \nonumber \\ 
&\qquad \bigg\|\exp\bigg(\int_0^T \|u_0\|_{\mathbb L^2}\|\nabla u\|_{\mathbb L^2} dr\bigg) \bigg\|_{\mathbb L^8(\Omega)} \nonumber \\
&\qquad\quad \bigg\|\exp\bigg(\int_0^T \|u_0^h\|_h\|\delta_+ u^h\|_h dr \bigg)\bigg\|_{\mathbb L^8(\Omega)}.
\end{align}
We proceed with bounding the three expectations on the right hand side of the above inequality. 

By Taylor expansion, there exists $\theta_l \in [(l-1)h,(l+1)h]$ such that
$R^h(l):=u_{xxxx}(\theta_l) h^2/4!$. 
Combining the Sobolev embedding \eqref{sob}, Corollary \ref{us-sup} implies, with $p=4$ and $\bs=5$, that
\begin{align}\label{e1}
\sqrt[4]{\ee\bigg[\sup_{t\in [0,T]} \|R^h\|_{l^\infty_h} ^4\bigg]}
\le C\sqrt[4]{\ee\bigg[\sup_{t\in[0,T]}\|u\|_{\hh^5}^4\bigg]} h^2
\le Ch^2.
\end{align}
By the exponential integrability \eqref{exp-mom-u} and \eqref{exp-mom-uh} in Theorem \ref{exp-u-uh}, we get
\begin{align}\label{e2}
&\bigg\|\exp\bigg(\int_0^T \|u_0\|_{\mathbb L^2}\|\nabla u\|_{\mathbb L^2} dr\bigg) \bigg\|_{\mathbb L^8(\Omega)} 
\bigg\|\exp\bigg(\int_0^T \|u_0^h\|_h\|\delta_+ u^h\|_h dr \bigg)\bigg\|_{\mathbb L^8(\Omega)}
<\infty.
\end{align}
We conclude \eqref{u-uh0} by combining \eqref{err-str}--\eqref{e2}.
\qed\\

A slight version of the proof of Theorem \ref{u-uh} leads to the following error estimate of $\epsilon^h$ under the $\mathbb L^p(\Omega; \CC([0,T]; \mathbb L^2))$-norm for all $p\ge 4$.

\begin{cor}\label{u-uhp}
Let $p=2$ or $p\ge 4$.
Assume that $u_0\in \bigcap\limits_{m=2}^5 \mathbb L^{2p\cdot 3^{5-m}}(\Omega; \hh^m)$, \eqref{con-u1}--\eqref{con-u2} hold for $q=4p$  and $Q^\frac12\in \LL_2^{5}$.
Then there exists a constant $C=C(T,p,u_0,Q)$ such that
\begin{align}\label{u-uh1}
\bigg(\ee\bigg[\sup_{t\in [0,T]}\|u(t)-u^h(t)\|_h^p\bigg]\bigg)^\frac1p
\le C h^2.
\end{align}
\end{cor}

\textbf{Proof.} 
Similarly to the proof of Theorem \ref{u-uh}, we apply the It\^o formula to $\|\epsilon(t)\|_h^p$, combined with the Young inequality, and obtain
\begin{align*}
& \|\epsilon(t)\|_h^p \\
&\le p\int_0^t \|\epsilon(r)\|_h^p
\Big(\|u_0\|_{\mathbb L^2}\|\nabla u\|_{H}+\|u_0^h\|_h\|\delta_+u^h\|_h\Big) dr+ p \int_0^t  ||\epsilon(r)||_h^{p-1}\|R^h\|_{l_h^\infty}dr \\
&\le \int_0^t \|R^h\|_{l_h^\infty}^pdr
+\int_0^t\|\epsilon(r)\|_h^p \Big(p-1+p\|u_0\|_{\mathbb L^2}\|\nabla u\|_{H}+p\|u_0^h\|_h\|\delta_+u^h\|_h\Big)dr.
\end{align*}
Applying Gronwall-Bellman inequality, and then taking expectation and using H\"older and Minkovskii inequalities, as in the proof of Theorem \ref{u-uh}, we deduce
\begin{align}\label{err-str}
\bigg(\ee\bigg[\sup_{t\in [0,T]}\|\epsilon^h(t)\|_h^p\bigg]\bigg)^\frac1p
&\le C \bigg(\ee\bigg[\sup_{t\in [0,T]} \|R^h\|_{l_h^\infty}^{2p}\bigg]\bigg)^\frac1{2p} \nonumber \\ 
&\qquad \bigg\|\exp\bigg(\int_0^T \|u_0\|_{\mathbb L^2}\|\nabla u\|_{\mathbb L^2} dr\bigg) \bigg\|_{\mathbb L^{4p}(\Omega)} \nonumber \\
&\qquad\quad \bigg\|\exp\bigg(\int_0^T \|u_0^h\|_h\|\delta_+ u^h\|_h dr \bigg)\bigg\|_{\mathbb L^{4p}(\Omega)}.
\end{align}
We conclude \eqref{u-uh1} by Corollary \ref{us-sup} and Theorem \ref{exp-u-uh}.
\qed\\

If $u_0$ is a non-random element in $\hh^\bs$ with $\bs\ge 2$, then the conditions \eqref{u0} and \eqref{con-u1}--\eqref{con-u2} disappear and we always obtain the exponential integrability theorem \eqref{exp-u-uh}.
In this case, our main results, Theorem \ref{u-uh} and Corollary \ref{u-uhp}, reduce to the following result.

\begin{tm}
Let $p\ge 1$.
Assume that $u_0\in \hh^5$ and $Q^\frac12\in \LL_2^{5}$.
Then there exists a constant $C=C(T,p,u_0,Q)$ such that
\begin{align*}
\bigg(\ee\bigg[\sup_{t\in [0,T]}\|u(t)-u^h(t)\|_h^p\bigg]\bigg)^\frac1p
\le C h^2.
\end{align*}
\end{tm}

\begin{rk}
\begin{enumerate}
\item
Our error analysis is also available under rough regularity assumptions.
More precisely, for some $\delta>1/2$ and $0<\beta<1$, if $u\in H^{2+\beta+\delta}$ (the fractional Sobolev space) a.s., then $\|R_h\|_{l_h^\infty} \le C \|u\|_{H^{2+\beta+\delta}} h^\beta$, which implies that the strong convergence rate is $\mathcal O(h^\beta)$; 
if $u\in H^{3+\beta+\delta}$ a.s., then 
$\|R_h\|_{l_h^\infty} \le C \|u\|_{H^{3+\beta+\delta}}h^{1+\beta}$ which yields that  the strong convergence rate is $\mathcal O(h^{1+\beta})$.

\item
If one have a priori estimate of $u$ under the $H^\delta$ norm for some $\delta>0$, we can also obtain strong convergence rate for our scheme \eqref{dif-nls} under weak regularity assumptions. 
For example, once a priori estimate under the $H^{4+\delta}$-norm with $\delta>1/2$ is established, we can reduce the regularity condition $H^5$ in Theorem \ref{u-uh} to $H^{4+\delta}$.
When the regularity exceeds $H^{4+\delta}$ with $\delta>1/2$, the order of the scheme \eqref{dif-nls} can not be improved. 
In this case, to obtain higher order schemes one can consider other higher order finite difference schemes rather than the central difference scheme \eqref{dif-nls} or use the extrapolation acceleration skill (see e.g. \cite{GK10}).
\end{enumerate}
\end{rk}

\section*{Acknowledgments}

The authors gratefully thank the anonymous referees for valuable comments and suggestions in improving this paper.

\bibliography{bib}

\end{document}